\documentclass[10pt]{amsart}

\usepackage{amssymb,amsmath}
\usepackage{amsthm,mathabx}
\usepackage{qsymbols}

\newtheorem{thm}{Theorem}
\newtheorem{lem}[thm]{Lemma}
\newtheorem{prop}[thm]{Proposition}

\newtheorem{defi}[thm]{Definition}

\newtheorem{remark}[thm]{Remark}

\newcommand{\X}{\mathcal{X}}

\newcommand{\R}{\mathbb{R}}

\newcommand{\E}{\mathbb{E}}
\newcommand{\ep}{\varepsilon}

\usepackage{color}

\begin{document}
\title{A proof of the Caffarelli contraction theorem via entropic regularization}

\begin{abstract}
We give a new proof of the Caffarelli contraction theorem, which states that the Brenier optimal transport map sending the standard Gaussian measure onto a uniformly log-concave probability measure is Lipschitz. The proof combines a recent variational characterization of Lipschitz transport map by the second author and Juillet with a convexity property of optimizers in the dual formulation of the entropy-regularized optimal transport (or Schr\"odinger) problem. 
\end{abstract}

\author{Max Fathi}
\author{Nathael Gozlan}
\author{Maxime Prodhomme}
\date{\today}

\date{\today}

\address{MF : Institut de Math\'ematiques de Toulouse, Universit\'e Paul Sabatier
118, route de Narbonne F-31062 Toulouse Cedex 9}
\email{mfathi@phare.normalesup.org}
\address{NG : Universit\'e Paris Descartes, MAP5, UMR 8145, 45 rue des Saints Pères, 75270 Paris Cedex 06}
\email{natael.gozlan@parisdescartes.fr}
\address{MP : Institut de Math\'ematiques de Toulouse, Universit\'e Paul Sabatier
118, route de Narbonne F-31062 Toulouse Cedex 9}
\email{maxime.prodhomme@math.univ-toulouse.fr}


\maketitle

\section{Introduction}

The aim of the paper is to give a new proof of the celebrated Caffarelli contraction theorem \cite{Caf00, Caf02}, which states that the Brenier optimal transport map sending the standard Gaussian measure on $\R^d$, denoted by $\gamma_d$ in all the paper, onto a probability measure $\nu$ having a log-concave density with respect to $\gamma_d$ is a contraction.
More precisely, let us recall the generalized version of Caffarelli's theorem:
\begin{thm}[Caffarelli \cite{Caf00, Caf02}]\label{thmCCT}For any probability measures $\mu$ and $\nu$ respectively of the form $\mu(dx) = e^{V(x)}\gamma_d(dx)$ and $\nu(dx) = e^{-W(x)}\gamma_d(dx)$ with $V$ and $W$ convex functions, and further assuming $\mu$ has a finite second moment and $\nu$ is compactly supported, there exists a continuously differentiable and convex function $\phi:\R^d \to \R$ such that $\nabla \phi$ is $1$-Lipschitz and $\nu = \nabla \phi_\# \mu.$
\end{thm}

Caffarelli's original result was only stated for the important particular case where $\mu$ is the Gaussian measure $\gamma_d$ (i.e. $V = 0$), but his proof readily extends to this more general setting \cite{Kol13}. Note that the assumption that $\nu$ is compactly supported can be removed via approximation. See \cite[Corollary 5.21]{Vil09} for details. Note that in all the paper we allow convex function to take the value $+\infty$.

This result plays an important role in the Functional Inequality literature, as it enables to transfer geometric inequalities such as Log-Sobolev or Gaussian Isoperimetric inequalities from the Gaussian measure to probability measures with a uniformly log-concave density. See \cite{Cor02,Har99, HAR01, Mil18} for some applications of Theorem \ref{thmCCT} to functional inequalities. It has also been used to derive deficit estimates in functional inequalities \cite{DPF17, CFP18}. Crucial for such applications is the dimension-free nature of the bound, to preserve the dimension-independent estimates that arise from these functional inequalities, and which are at the center of their applications in statistics for example. More recently, there have been some extensions, such as dimension-free Sobolev estimates \cite{Kol10, Kol13} and variants for compactly-supported perturbations of the Gaussian measure \cite{CFJ17}. 

Caffarelli's original proof relies on the formulation of Brenier maps as solutions to a Monge-Amp\`ere equation, and uses maximum principle-type estimates. In particular, it does not actually exploit the fact that $\nabla \phi$ is an optimal transport map. This is also the case for the other proofs \cite{Kol13, KM12}. Our purpose here is to provide a different proof that does directly exploit ideas from optimal transport theory. 

In this paper, we develop an approach based on a variational characterization of Lipschitz regularity of optimal transport maps obtained by the second author and Juillet in \cite{GJ18}. To recall this result, we need to introduce some notations and definitions.
We will denote by $\mathcal{P}(\R^d)$ the set of Borel probability measures on $\R^d$ and by $\mathcal{P}_{k} (\R^d)$, $k\geq1$, the subset of $\mathcal{P}(\R^d)$ of probability measures having a finite moment of order $k$. The quadratic Kantorovich distance $W_2$ is defined for all $\mu,\nu \in \mathcal{P}_2(\R^d)$ as follows:
\[
W_2^2 (\mu,\nu) = \inf_{\pi \in C(\mu,\nu)} \int |x-y|^2\,\pi(dxdy),
\]
where $|\,\cdot\,|$ denotes in all the paper the standard Euclidean norm and $C(\mu,\nu)$ is the set of couplings between $\mu$ and $\nu$, that is to say the set of probability measures on $\R^d\times \R^d$ such that $\pi(A\times \R^d) = \mu(A)$ and $\pi(\R^d \times B) = \nu(B)$, for all Borel sets $A,B$  of $\R^d.$
Finally, if $\eta_1,\eta_2 \in \mathcal{P}_1(\R^d)$, one says that $\eta_1$ is dominated by $\eta_2$ for the convex order if $\int f\,d\eta_1 \leq \int f\,d\eta_2$ for all convex function $f:\R^d \to \R$. In this case, we write $\eta_1 \leq_c \eta_2.$
With these notations in hand, the variational characterization of \cite{GJ18} reads as follows:
\begin{thm}\label{thmGJ}
Let $\mu,\nu \in \mathcal{P}_2(\R^d)$ ; the following assertions are equivalent:
\begin{itemize}
\item[(i)] There exists a continuously differentiable and convex function $\phi:\R^d \to \R$ such that $\nabla \phi$ is $1$-Lipschitz and $\nu = \nabla \phi_\# \mu$,
\item[(ii)] For all $\eta \in \mathcal{P}_2(\R^d)$ such that $\eta \leq_c \nu$, 
\[
W_2(\nu,\mu) \leq W_2(\eta,\mu).
\]
\end{itemize}
\end{thm}
In other words the Brenier map between $\mu$ and $\nu$ is a contraction if and only if $\nu$ is the closest point to $\mu$ among all probability measures dominated by $\nu$ in the convex order.

Our strategy to recover Theorem \ref{thmCCT} is thus to show the following monotonicity property of the $W_2$ distance: if $\mu$ and $\nu$ satisfy the assumptions of Theorem \ref{thmCCT}, it holds
\begin{equation}\label{eq:minW_2}
W_2 (\nu,\mu) \leq W_2(\eta,\mu),\qquad \forall \eta \leq_c \nu.
\end{equation}
For that purpose, we will establish a similar inequality at an information theoretic level replacing $W_2$ by the so called entropic transport cost $\mathcal{T}_H^\ep$ (presented in details in the next section) that is defined in terms of the minimization of the relative entropy between $\pi$ and a reference measure $R^\ep$ involving some small noise parameter $\ep>0$. We will prove the following monotonicity property of the entropic cost:
\begin{equation}\label{eq:minT_H}
\mathcal{T}_H^\ep(\nu,\mu) \leq \mathcal{T}_H^\ep(\eta,\mu)
\end{equation}
for all $\eta \leq _c \nu$ with a finite Shannon information.
As observed by Mikami \cite{Mik04} and extensively developed by L\'eonard \cite{Leo12,Leo14} the zero noise limit of $\ep\mathcal{T}_H^\ep$ is $\frac{1}{2}W_2^2$. Thus letting $\ep \to 0$ in \eqref{eq:minT_H} will give \eqref{eq:minW_2}.

The paper is organized as follows. Section 2 introduces entropic transport costs and the main results of the paper. Section 3 gives proofs of these results.

\section{Entropic transport costs and main results}
\subsection{Entropic costs and their zero-noise limit}
Consider the classical Ornstein-Uhlenbeck process $(Z_t)_{t\geq0}$ on $\R^d$, defined by the following stochastic differential equation:
\[
dZ_t = \frac{-1}{2}Z_t\,dt + \,dW_t,\qquad t\geq0,
\]
where $(W_t)_{t\geq0}$ is a standard $d$ dimensional Brownian motion and $Z_0^\ep \sim \gamma_d$.
As it is well known, the process $Z$ admits the following explicit representation
\[
Z_t = Z_0e^{- t/2} + e^{- t/ 2} \int_0^t e^{ s /2}\,dW_s,\qquad t\geq0.
\]
The joint law of $(Z_0, Z_\ep)$ will be denoted by $R^\ep$.
It is therefore given by
\[
R^\ep = \mathrm{Law} \left(X , X e^{-\ep /2} + \sqrt{1 - e^{-\ep  } }Y\right),
\]
with $X,Y$ two independent standard Gaussian random vectors on $\R^d$. In other words, 
\[
R^\ep(dxdy)  = \gamma_d(dx) r_x^\ep(dy),
\]
where $x\mapsto r_x^\ep$ is the probability kernel defined by $r_x^\ep = \mathcal{N}(x e^{-\ep  /2}, (1-e^{-\ep }) I_d)$.

Recall that the relative entropy of a probability measure $\alpha$ with respect to another probability measure $\beta$ on some measurable space $(\X,\mathcal{A})$ is defined by
\[
H(\alpha | \beta) = \int \log \left(\frac{d\alpha}{d\beta}\right)\,d\alpha,
\]
if $\alpha$ is absolutely continuous with respect to $\beta$. If this is not the case, one sets $H(\alpha|\beta)=+\infty.$ The relative entropy is a non-negative quantity that vanishes only when the two probability measures are equal, this is why it is often called Kullback-Leibler distance (even though it is not a true distance).
\begin{defi}[Entropic transport cost]
For all probability measures $\mu,\nu$ on $\R^d$, the entropic transport cost associated to $R^\ep$ is defined by
\[
\mathcal{T}_H^\ep(\mu,\nu) = \inf_{\pi \in \Pi(\mu,\nu)} H(\pi| R^\ep).
\]
\end{defi}
As shown by Mikami, L\'eonard and others \cite{Mik04,Leo12} the zero noise limit of $\ep \mathcal{T}_H^\ep$ is $\frac{1}{2}W_2^2$. At a heuristic level, this phenomenon can be easily understood from the following identities:
\begin{align*}
\ep H(\pi | R^\ep) &= \ep \int \log\left(\frac{d\pi}{dx}\right)\,d\pi - \ep \int \log\left(\frac{dR^\ep}{dx}\right)\,d\pi\\
&= \ep\int \log\left(\frac{d\pi}{dx}\right)\,d\pi  + \frac{\ep}{2(1-e^{-\ep})} \int |y- e^{-\ep/2}x|^2\,\pi(dxdy) + \frac{\ep}{2} \int |x|^2\,\mu(dx) + c(\ep), 
\end{align*}
where $c(\ep) \to 0$ (and is independent of $\mu,\nu,\pi$). So for small $\ep$, minimizing $\pi \mapsto H(\pi|R^\ep)$ amounts to minimizing $\pi \mapsto \frac{1}{2}\int |x-y|^2\,\pi(dxdy).$  

In the sequel we will use the following result, which can be easily deduced from a general convergence theorem due to Carlier, Duval, Peyr\'e and Schmitzer \cite[Theorem 2.7]{CDPS17}.
We will say that a probability measure $\eta$ is of finite (Shannon) entropy if it is absolutely continuous with respect to the Lebesgue measure and if $\int \log \left(\frac{d\eta}{dx}\right)\,d\eta$ is finite. Note that, if $\eta \in \mathcal{P}_2(\R^d)$, then it is of finite entropy if and only if $H(\eta|\gamma_d)<+\infty.$

\begin{thm}[Carlier et al. \cite{CDPS17}]\label{thm:limit}
Suppose that $\mu,\nu \in \mathcal{P}_2(\R^d)$ are of finite entropy. Then, it holds
\[
\ep \mathcal{T}_H^\ep(\mu,\nu) \to \frac{1}{2}W_2^2(\mu,\nu)\qquad \text{as } \ep \to 0.
\]
\end{thm}
We state now a technical lemma that will be needed to apply Theorem \ref{thm:limit} in our framework:
\begin{lem}\label{finite-entropy}If $\mu$ and $\nu$ satisfy the assumptions of Theorem \ref{thmCCT}, then they are of finite entropy.
\end{lem}
The proof is postponed to Section \ref{Sec:Proofs}.

\subsection{Entropic cost in the framework of Caffarelli theorem}
As explained above, the key step in our approach consists in showing that on the set of probability measures dominated by $\nu$ in the convex order, the closest point to $\mu$ for the entropic cost distance is $\nu$ itself (when $\nu$ satisfies the assumptions of Theorem \ref{thmCCT}). 

\begin{thm}\label{thm:minT_H}
Let $\mu$ and $\nu$ satisfy the assumptions of Theorem \ref{thmCCT}. Additionally assume that $V$ is bounded from below. If $\eta$ is such that $\eta \leq_c  \nu$, then for all $\ep>0$
\[
\mathcal{T}_H^\ep (\mu,\eta ) \geq \mathcal{T}_H^\ep (\mu,\nu).
\]
\end{thm}
Let us admit Theorem \ref{thm:minT_H} for a moment and complete the proof of Theorem \ref{thmCCT}.
\proof[Proof of Theorem \ref{thmCCT}]

Let us temporarily assume that $V$ is bounded from below. According to Lemma \ref{finite-entropy}, $\mu$ and $\nu$ have finite entropy. So using Theorem \ref{thm:limit}, one concludes by letting $\ep \to 0$ that for all compactly supported probability measures $\nu$ of the form $\nu(dx) = e^{-W(x)}\,\gamma_d(dx)$, with $W:\R^d\to \R\cup\{+\infty\}$ convex, it holds
\[
W_2(\mu,\nu) \leq W_2(\mu,\eta)
\] 
for all $\eta$ of finite entropy and such that $\eta\leq_c\nu$. Now, fix some compactly supported $\nu_0$ of the form $\nu_0(dx) = e^{-W_0(x)}\,\gamma_d(dx)$, with $W_0:\R^d\to \R\cup\{+\infty\}$ convex and let us show that the inequality \eqref{eq:minW_2} holds for any $\eta \leq_c \nu$. Take $\eta \leq_c \nu_0$ and define, for all $\theta \in (0,\pi/2)$,  
\[
\nu_\theta = \mathrm{Law} (\cos \theta X + \sin \theta Z)\qquad \text{and} \qquad \eta_\theta = \mathrm{Law} (\cos \theta Y + \sin \theta Z),
\]
where $X\sim \nu_0$, $Y \sim \eta$ and $Z$ is independent of $X$ and $Y$ and has density $\frac{1}{C} \mathbf{1}_B(x) e^{-\frac{|x|^2}{2}}$, where $B$ is the Euclidean unit ball. According to Lemma \ref{lem:technique}, $\nu_\theta$ is compactly supported and of the form $e^{-W_\theta}\gamma_d$, with $W_\theta$ convex, $\eta_\theta$ is of finite entropy and $\eta_\theta \leq_c \nu_\theta$. Therefore, it holds $W_2(\mu,\nu_\theta) \leq W_2(\mu,\eta_\theta)$. Letting $\theta \to 0$ gives $W_2(\mu,\nu_0) \leq W_2(\mu,\eta)$ for all $\eta \leq_c \nu_0$, which, according to Theorem \ref{thmGJ}, completes the proof when $\mu$ has finite entropy. 

Finally, let us remove the assumption that $V$ is bounded from below. 
%
Since $V$ is convex, it is bounded from below by some affine function. Thus there exists some $a \in \R^d$ such that $x\mapsto V(x) + a \cdot x$ is bounded from below. Consider the probability measure $\tilde{\mu}$ defined as the push forward of $\mu$ under the translation $x \mapsto x+a$. An easy calculation shows that the density of $\tilde{\mu}$ with respect to $\gamma_d$ is $Ce^{V(x-a) + a\cdot (x-a)}$, with $C$ a normalizing constant, and so $\tilde{\mu}$ satisfies our assumptions. Therefore, there exists a continuously differentiable convex function $\tilde{\phi}:\R^d \to \R$ such that $\nabla \tilde{\phi}$ is $1$-Lipschitz and $\nu = \nabla\tilde{\phi}_\# \tilde{\mu}$. Setting $\phi(x) = \tilde{\phi}(x+a)$, $x \in \R^d$, one gets $\nu = \nabla\phi_\# \mu$ which completes the proof.
\endproof

Before proving Theorem \ref{thm:minT_H}, let us informally explain why one can guess the statement is easier to prove at the level of entropic cost than directly for the Wasserstein distance. If we consider the plain relative entropy, we have the variational formula 
\[
H(\rho | \mu) = \sup \int{f \,d\rho} - \log \int{e^f \,d\mu},
\]
where the supremum runs over the set of functions $f$ such that $\int e^{f}\,d\mu<+\infty$.
Hence, taking $f = -(V+W)$, gives
\[
H(\rho | \mu) \geq \int{-(V + W) \,d\rho} \geq \int{-(V + W)\, d\nu} = H(\nu | \mu)
\]
as soon as $\rho \leq_c \nu$. So this trivial bound hints at the fact that comparison is easier for entropies when we have a log concavity condition on the relative density. 

To prove Theorem \ref{thm:minT_H}, we need to know more about the optimal coupling $\pi$ for $\mathcal{T}_H^\ep (\mu,\nu)$. 
The following is classical in entropic transport literature and goes back to the study of the so called Schr\"odinger bridges \cite{Sch32}.
\begin{prop}\label{prop:piopt}
Let $\mu,\nu \in \mathcal{P}_2(\R^d)$ be such that $H(\mu|\gamma_d)<+\infty$ and $H(\nu|\gamma_d)<+\infty$
\begin{enumerate}
\item There exists 
a unique coupling $\pi^\ep \in C(\mu,\nu)$ such that 
\[
\mathcal{T}_H^\ep(\mu,\nu) = H(\pi^\ep | R^\ep) <+\infty
\]
\item There exist two measurable functions $f^\ep,g^\ep : \R^d \to \R^+$ such that $\log f^\ep \in L^1(\mu)$, $\log g^\ep \in L^1(\nu)$ and
\[
\pi^\ep(dxdy) = f^\ep(x)g^\ep(y)R^\ep(dxdy).
\] 
\end{enumerate}
\end{prop}
\proof[Sketch of proof] (1) We equip the set $\mathcal{P}(\R^d\times \R^d)$ with the usual topology of narrow convergence. For this topology, the function $\pi \mapsto H(\pi |R^\ep)$ is lower-semicontinuous and the set $C(\mu,\nu)$ is compact. Therefore, the function $H(\,\cdot\,|R^\ep)$ attains its minimum at some point $\pi^\ep$ of $C(\mu,\nu)$. It is easily checked that the coupling $\pi_0 = \mu \otimes \nu$ is such that $H(\pi_0 | R^\ep) <+\infty$, so $H(\pi^\ep | R^\ep) <+\infty$. Uniqueness comes from the strict convexity of $H(\,\cdot\,|R^\ep)$. For the proof of (2) we refer to \cite[Corollary 3.2]{Csi75}. In the special case where $\mu$ and $\nu$ satisfy our log-convexity/concavity assumptions we will give a self-contained proof in Section \ref{Sec:Proofs}.
\endproof

%

In the setting of Theorem \ref{thmCCT}, it turns out that much more can be said about the functions $f$ and $g$. This is explained in the following result, which seems of independent interest.

\begin{thm}\label{thm:minimizer} With the same notation as in Proposition \ref{prop:piopt}, let $\mu$ be a probability measure of the form $\mu(dx)=e^{V(x)}\gamma_d(dx)$ with a finite second moment and $\nu$ be a compactly supported probability measure on $\R^d$ of the form $\nu(dx) = e^{-W(x)}\,\gamma_d(dx)$, with $V, W$ convex and $V$ bounded from below. There exist a \emph{log-convex} function $f^\ep:\R^d \to [1,+\infty)$ and a \emph{log-concave} function $g^\ep: \R^d \to [0,+\infty)$ such that the unique optimal coupling $\pi^\ep \in \Pi(\mu,\nu)$  is of the form $\pi^\ep(dxdy) = f^\ep(x)g^\ep(y)\,R^\ep(dxdy)$. Moreover, the function $\log f^\ep$ is integrable with respect to $\mu$ and the function $\log g^\ep$ is integrable with respect to $\nu$ and it holds
\[
\mathcal{T}_H^\ep(\mu,\nu) = H(\pi^\ep | R^\ep) = \int \log f^\ep \,d\mu + \int \log g^\ep \,d\nu.
\]
\end{thm}


We now give a brief heuristic explanation as to why one can expect this statement to imply the Caffarelli contraction theorem. Informally, from the convergence of the entropic cost to the Wasserstein distance, we expect from the dual formulation that $\epsilon \log f$ converges to $|x|^2/2 - \varphi$ (up to some additive constant), where $\varphi$ is a potential giving rise to the optimal transport map $T = \nabla \varphi$. However convexity is preserved by pointwise convergence, so we expect $|x|^2/2 - \varphi$ to also be convex. But this is equivalent to $\nabla \varphi$ being 1-Lipschitz, since the eigenvalues of the Hessian of $\varphi$ must then be bounded by 1. Theorem \ref{thmGJ} will allow us to avoid having to prove convergence of $\epsilon \log f$ to a Kantorovitch potential. 

Section \ref{Sec:Proofs} is essentially devoted to the proof of Theorem \ref{thm:minimizer}. With Theorem \ref{thm:minimizer} in hand, the proof of Theorem \ref{thm:minT_H} becomes almost straightforward:
\proof[Proof of Theorem \ref{thm:minT_H}] 
Recall the following duality inequality for the relative entropy : if $\alpha,\beta$ are two probability measures on a measurable space $(\mathcal{X},\mathcal{A})$ such that $H(\alpha | \beta)<+\infty$, then for any measurable function $h: \X \to \R$ such that  $\int e^h\,d\beta <+\infty$, it holds $\int [h]_+ d\alpha <+\infty$ and
\begin{equation}\label{eq:dual}
H(\alpha | \beta) \geq \int h\,d\alpha - \log \left( \int e^h\,d\beta\right)
\end{equation}
Let $\pi \in C(\mu,\eta)$ be a coupling between $\mu$ and some probability $\eta \leq_c \nu$ such that $H(\pi | R^\ep)<+\infty$ ;  applying the inequality above to $\alpha = \pi$, $\beta = R^\ep$ and $h(x,y) = \log \left(f^\ep(x)g^\ep(y)\right)$, $x,y \in \R^d$ gives
\begin{align*}
H(\pi |R^\ep) &\geq \int \log f^\ep (x) + \log g^\ep (y) \,\pi(dxdy)\\
& = \int \log f^\ep (x)\,\mu(dx) + \int \log g^\ep (y) \,\eta(dy)\\
& \geq \int \log f^\ep (x)\,\mu(dx) + \int \log g^\ep (y) \,\nu(dy)\\
& = H(\pi^\ep | R^\ep) = \mathcal{T}^\ep_H(\mu,\nu),
\end{align*}
where the second inequality comes from the fact that $\log g^\ep$ is a \emph{concave} function and $\eta \leq_c \nu.$ Optimizing over $\pi$, gives the inequality $\mathcal{T}^\ep_H(\mu,\eta) \geq \mathcal{T}^\ep_H(\mu,\nu)$ and completes the proof.
\endproof

To conclude this section, we mention some perspectives. The most natural question is whether this scheme of proof can be adapted to establish a version of Caffarelli's theorem in other settings than $\R^d$, such as on manifolds or in free probability \cite{GS14}. See \cite{Mil18} for some motivations in analysis and geometry. See \cite{GT18} for a study of Schr\"odinger's problem in a wider geometric setting. Another question is about integrated or non-local quantitative regularity estimates, such as those in \cite{Kol10, Kol13}. The role of 1-Lipschitz bounds in Theorem \ref{thmGJ} is very specific, we do not know if there is an analogue of that equivalence adapted to other types of regularity bounds. However, it could be possible to prove stable a priori bounds on $\ep \log f^\ep$ and pass to the limit. Of particular interest is whether we can establish integrated gradient bounds for non-uniformly convex potentials, since such estimates can still be used to establish Poincar\'e inequalities \cite{Mil09,Kla13}. Finally, \cite{DPF17} proves a rigidity/stability result for the Caffarelli contraction theorem, and it would be interesting to find a way to improve the quantitative bounds.

\section{Proofs} \label{Sec:Proofs}
This section contains the material needed to prove Theorem \ref{thm:minT_H}. The ideas developed here are adapted from a paper by Fortet \cite{For40}. We warmly thank Christian L\'eonard for mentioning us this paper and explaining to us the ingredients of Fortet's proof. Fortet's work was also recently revisited in \cite{EP19}. 

We will denote by $P^\ep$ the Ornstein-Uhlenbeck semi-group at time $\ep$ defined for all non-negative measurable function $\psi$ by
\[
P^\ep \psi (x) = \E[\psi (Z_\ep) | Z_0=x] = \frac{1}{(2\pi)^{d/2}} \frac{1}{(1-e^{-\ep })^{d/2}} \int_{\R^d} \psi(y + e^{-\ep  /2} x) e^{- \frac{|y|^2}{2(1-e^{-\ep })}} \,dy,\qquad x \in \R^d.
\]

Suppose that $f^\ep,g^\ep$ are measurable non-negative functions such that $\pi^\ep (dxdy)= f^\ep(x)g^\ep(y)R^\ep(dxdy)$ belongs to $C(\mu,\nu)$. Then, writing the marginals condition, one sees that $f^\ep$ and $g^\ep$ are related to each other by the identities: for all $x,y \in \R^d$
\begin{equation}\label{eq:f,g}
f^\ep(x) P^{\ep}g^\ep(x) = e^{V(x)}\qquad \text{and} \qquad g^\ep(y)P^\ep f^\ep(y) = e^{-W(y)}.
\end{equation}

These relations suggest to introduce the functional $\Phi^\ep$ defined as follows: for all measurable function $h: \R^d \to \R\cup \{+\infty\}$, 
\[
\Phi^\ep(h) = V -\log \left(P^{\ep} \left(e^{-W} \frac{1}{P^{\ep} (e^h)}\right)\right).
\]
With this notation, a couple $(f^\ep,g^\ep)$ satisfies \eqref{eq:f,g} if and only if $g^{\ep} = e^{-W} \frac{1}{P^\ep(f^\ep)}$ and $f^\ep = e^{h^\ep}$ with $h^\ep$ such that
\[
h^\ep = \Phi^\ep(h^{\ep}).
\]
This fixed point equation suggests that the unknown function $h^\ep$ could be obtained as the limit when $n \to +\infty$ of a sequence $(h_n)_{n\geq0}$ satisfying the recursive scheme 
\begin{equation}\label{eq:rec}
h_{n+1} = \Phi^\ep(h_n),\qquad n\geq0
\end{equation} 
and initialized with some function $h_0.$ This fixed point scheme is actually at the core of the use of Sinkhorn's algorithm to numerically approximate optimal transport via entropic regularization \cite{Cut13, BCC15}. 

The convexity of $h^\ep$ can then be established if we can initiate this fixed point scheme \eqref{eq:rec} with some convex initial data $h_0$, thanks to the following key result:
\begin{lem}\label{lem:key} If $h : \R^d \to \R\cup\{+\infty\}$ is convex, then $\Phi^\ep(h)$ is also convex.
\end{lem}
\proof
This property is inherited from the following classical properties of $P^\ep$ :
\begin{itemize}
\item If $f$ is log-convex, then $P^\ep(f)$ is log-convex. This simply follows from H\"older inequality.
\item If $g$ is log-concave, then  $P^\ep(g)$ is log-concave. This follows from the fact that the set of log-concave functions is stable under convolution which is a well known consequence of Prekopa Theorem \cite{Pre73}.
\end{itemize}
\endproof

The line of reasoning sketched above is essentially the one adopted in the proof of Theorem \ref{thm:minimizer}, except that the recurrence scheme  \eqref{eq:rec} needs to be properly modified in order to force its convergence (this modification is the same as the one proposed by Fortet in \cite{For40}).

\begin{remark}
In the compact setting, the map $\Phi^\ep$ is actually a contraction with respect to a well-chosen metric, see for example \cite[Lemma 1]{Gen19} or \cite{CGP16} (following the earlier ideas of \cite{FL89} in the discrete setting). This would ensure that the fixed point must belong to any stable subspace. Here, we work in a noncompact setting ($\mu$ has non-compact support) and it seems the map is globally only 1-Lipschitz at that level of generality. One could however expect that it remains a contraction on a suitable stable subspace of convex functions. 
\end{remark}

\begin{remark}
A natural question is whether our scheme of proof can be used directly at the level of the Kantorovitch dual formulation of optimal transport, rather than on the regularized version. The answer seems to be no, as in the limit while the minimizers in the dual formulation of entropic transport, suitably rescaled, converge to the Kantorovitch potentials, the fixed point problem becomes degenerate in the limit, and only selects so-called c-convex functions (with the cost here being the quadratic distance), so we lose uniqueness. Indeed, there is no known fixed point scheme similar to \eqref{eq:rec} for Kantorovitch potentials, which is why Sinkhorn's algorithm is only used to numerically approximate the regularized problem \cite{Cut13}. 
\end{remark}

Before moving on to the proof, let us present two other essential properties of $\Phi^\ep$.
\begin{lem}\label{lem:propPhi}\ 
\begin{enumerate}
\item The map $\Phi^\ep$ is monotone: $h\leq k \Rightarrow \Phi^{\ep} (h) \leq \Phi^\ep(k).$
\item For any measurable $h:\R^d \to \R$, it holds
\[
\int \exp\left(h(x)- \Phi^\ep(h)(x)\right) d\mu \leq 1,
\]
with equality if $h$ is bounded from above.
\end{enumerate}
\end{lem}
\proof
The first point is straightforward. Let us prove the second point. Since the operator $P^\ep$ is symmetric in $L^2(\gamma_d)$, for any function $h:\R^d \to \R$ it holds
\[
\int e^{h\wedge a  - \Phi^\ep(h)}\,d\mu = \int e^{h\wedge a} P^\ep\left(e^{-W} \frac{1}{P^\ep(e^{h})}\right)\,d\gamma_d
 = \int P^\ep\left(e^{h\wedge a}\right) e^{-W} \frac{1}{P^\ep(e^{h})}\,d\gamma_d.
\]
Letting $a \to +\infty$, one gets by monotone convergence $\int e^{h - \Phi^\ep(h)}\,d\mu = \int_{\{P^\ep(e^h)<+\infty\}} e^{-W}\,d\gamma_d$ which gives the claim.
\endproof

The existence of a coupling of the desired form can be established under more general conditions on $\mu$ and $\nu$:
\begin{thm}\label{thm:existence} Let $\mu$ be a probability measure of the form $\mu(dx) = e^{V(x)}\gamma_d(dx)$ with $V:\R^d \to \R$ convex and bounded from below, and let $\nu$ be a probability measure on $\R^d$ of the form $\nu(dx) = e^{-W(x)}\,\gamma_d(dx)$, with $W : \R^d \to \R\cup\{+\infty\}$ a convex function such that $\{W< -m\}$ is bounded for $m = \inf_{\R^d} V \leq0$. There exist a \emph{log-convex} function $f^\ep:\R^d \to [1,+\infty)$ and a \emph{log-concave} function $g^\ep: \R^d \to [0,+\infty)$ such that the measure $\pi^\ep$ defined by $\pi^\ep(dxdy) = f^\ep(x)g^\ep(y)\,R^\ep(dxdy)$ belongs to $C(\mu,\nu).$ 
\end{thm}

\proof[Proof of Theorem \ref{thm:existence}] 
Let us show that there exists a convex function $\bar h : \R^d \to \R^+$ such that $\Phi^\ep(\bar h) = \bar h$. Then, defining $f^\ep = e^{\bar h}$ and $g^\ep= e^{-V}/P^\ep(f^\ep)$, we see that $f^\ep$ is log-convex, $g^\ep$ is log-concave (we use again the fact that $P^\ep$ preserves log-convexity) and satisfy \eqref{eq:f,g}.

Let us define by induction the sequence $(h_n)_{n\geq0}$ as follows: $h_0=0$ and for all $n\geq0$
\begin{equation}\label{eq:suite}
h_{n+1} = \left[\Phi^\ep(h_n)\right]_+ \wedge n.
\end{equation}
By construction, note that $h_n$ takes values in $[0,n-1]$.
Let us show by induction that the sequence $(h_n)_{n\geq0} $ is non-decreasing. First observe that $h_1 = 0= h_0$ and so in particular $h_0\leq h_1$. According to Item (1) of Lemma \ref{lem:propPhi}, the operator $\Phi^\ep$ is non-decreasing. Therefore, if $h_{n+1} \geq h_n$ for some $n\geq0$, then 
\[
h_{n+2} = \left[\Phi^\ep(h_{n+1})\right]_+ \wedge (n+1) \geq \left[\Phi^\ep(h_n)\right]_+ \wedge (n+1) \geq \left[\Phi^\ep(h_n)\right]_+ \wedge n =h_{n+1}. 
\]
Let us denote by $h_\infty$ the pointwise limit of $h_n$ as $n \to \infty.$ The function $h_\infty$ takes values in $\R^+\cup \{+\infty\}.$
Let us show that $h_\infty$ solves the following fixed point equation
\begin{equation}\label{eq:fixed-point}
h_\infty = \left[\Phi^\ep(h_\infty)\right]_+.
\end{equation}
Indeed, by monotone convergence, $P^\ep(e^{h_n}) \to P^\ep(e^{h_\infty})$. Then,  by dominated convergence,
\[
P^\ep\left(e^{-W}\frac{1}{P^\ep(e^{h_n})} \right) \to P^\ep\left(e^{-W}\frac{1}{P^\ep(e^{h_\infty})} \right)
\]
which implies that $h_{n} \to  \left[\Phi^\ep(h_\infty)\right]_+$  and gives \eqref{eq:fixed-point}.

Now let us show that $h_\infty$ is in fact a fixed point of $\Phi^\ep$. Let us admit for now that $\Phi^\ep(h_\infty)(x)<+\infty$ for all $x \in \R^d.$ Then, thanks to \eqref{eq:fixed-point}, $h_\infty(x)<+\infty$ for all $x \in \R^d.$ According to Item (3) of Lemma \ref{lem:propPhi}, it holds $\int e^{h_\infty - \Phi^\ep(h_\infty)}\,d\mu \leq 1$. Since $h_\infty \geq \Phi^\ep(h_\infty)$, the function $e^{h_\infty - \Phi^\ep(h_\infty)}$ is bounded from below by $1$. Therefore, $h_\infty = \Phi^\ep(h_\infty)$ almost everywhere. The function $ \Phi^\ep(h_\infty)$ is easily seen to be continuous and since $h_\infty$ satisfies \eqref{eq:fixed-point} it is also continuous. The functions $h_\infty$ and $ \Phi^\ep(h_\infty)$ being continuous, the equality $h_\infty = \Phi^\ep(h_\infty)$ holds in fact everywhere. 
To complete the proof that $h_\infty$ is a fixed point of $\Phi^\ep$, it remains to prove that  $\Phi^\ep(h_\infty)(x)<+\infty$ for all $x \in \R^d.$ Let us assume, by contradiction, that there exists some $x_o\in \R^d$ such that $\Phi^\ep(h_\infty)(x_o)=+\infty.$ This easily implies that $P^\ep(e^{h_\infty}) = +\infty$ almost everywhere, which in turn implies that $\Phi^\ep(h_\infty)\equiv \infty$. Since $h_\infty \geq \Phi^\ep(h_\infty)$, one concludes also that $h_\infty \equiv +\infty.$ Now let us show that there exists $n_0$ such that for all $n \geq n_0$
\begin{equation}\label{eq:conv-unif}
\inf_{x \in \R^d} \Phi^\ep(h_n)(x)\geq0.
\end{equation}
For any $x \in \R^d$, it holds (denoting by $C= (2\pi)^{d/2}(1-e^{-\ep })^{d/2}$ and by $m = \inf_{\R^d} V$)
\begin{align*}
P^\ep\left(e^{-W} \frac{1}{P^\ep(e^{h_n})}\right)(x) &= \frac{1}{C} \int_{\R^d} e^{-W(y)} \frac{1}{P^\ep(e^{h_n})}(y) e^{- \frac{|y-e^{-\ep  /2} x|^2}{2(1-e^{-\ep })}} \,dy\\
& =  \frac{1}{C} \int_{\{W\leq -m\}} e^{-W(y)} \frac{1}{P^\ep(e^{h_n})}(y) e^{- \frac{|y-e^{-\ep  /2} x|^2}{2(1-e^{-\ep })}} \,dy\\
& + \frac{1}{C}  \int_{\{W> -m\}} e^{-W(y)} \frac{1}{P^\ep(e^{h_n})}(y) e^{- \frac{|y-e^{-\ep  /2} x|^2}{2(1-e^{-\ep })}} \,dy\\
& =  \frac{1}{C} \int_{\{W\leq -m\}}  e^{- \frac{|y-e^{-\ep  /2} x|^2}{2(1-e^{-\ep })}} \,dy \max_{z \in \{W\leq -m\}} e^{-W(z)} \frac{1}{P^\ep(e^{h_n})}(z).\\
& +\frac{e^{m}}{C}\int_{\{W> -m\}}  e^{- \frac{|y-e^{-\ep  /2} x|^2}{2(1-e^{-\ep })}} \,dy,
\end{align*}
where we used the fact that $P^\ep(e^{h_n})\geq 1$, since $h_n\geq0$. The sequence of functions $ \frac{1}{P^\ep(e^{h_n})}$ is a non-increasing sequence of continuous functions converging to $0$. Therefore, according to Dini's Theorem, the convergence is uniform on the compact set $K = \overline{\{W< 0\}}$. Since $W$ is convex, it is bounded from below on $K$. Therefore, there exists $n_0$ such that $\max_{z \in K} e^{-W(z)} \frac{1}{P^\ep(e^{h_n})}(z) \leq e^m$ for all $n\geq n_0.$ Plugging this inequality into the inequality above, one easily gets \eqref{eq:conv-unif}. Now, according to \eqref{eq:conv-unif}, there exists some $n_o$ such that $\Phi^\ep(h_{n_0}) \geq 0.$ Therefore $h_{n_o+1} = \Phi^\ep(h_{n_o}) \leq \Phi^\ep(h_{n_o+1}).$ Since $h_{n_0+1}$ is bounded, Item (2) of Lemma \ref{lem:propPhi} yields $\int e^{h_{n_o+1} - \Phi^\ep(h_{n_o+1})}d\mu=1$, which implies that $h_{n_o+1}= \Phi^\ep(h_{n_o+1})$. Therefore, $h_\infty = h_{n_o+1}$, which contradicts the fact that $h_\infty \equiv +\infty.$

Finally, let $k_0 = h_\infty^{**}$ be the convex regularization of $h_\infty$ (which is well defined since $h_\infty$ is bounded from below). By definition $ k_0 \leq h_\infty$ and since $h_\infty \geq0$, it holds $k_0 \geq 0$. 
Define by induction $(k_n)_{n \geq1}$ by $k_{n+1} = \max (\Phi^\ep(k_n) ; k_0)$. 
Since according to Lemma \ref{lem:key} $\Phi^\ep$ preserves convexity and $k_0$ is convex, $k_n$ is convex for all $n$.
The sequence $k_n$ is non-decreasing and satisfies $k_n \leq h_\infty$ for all $n$. Therefore, $k_n$ converges pointwise to some $k_\infty$, which is also convex and finite valued. Reasoning as above one sees that $k_\infty = \max (\Phi^\ep(k_\infty) ; k_0)$ and so in particular $k_\infty \geq \Phi^\ep(k_\infty).$ Using again the fact that $\int e^{k_\infty - \Phi^\ep(k_\infty)}\,d\mu \leq 1$, one concludes that $k_\infty$ is a fixed point of $\Phi^\ep$. Setting $\bar h = k_\infty$ completes the proof.
\endproof

\proof[Proof of Theorem \ref{thm:minimizer}]
First, let us note that $H(\mu\otimes \nu | R^\ep)<+\infty$. Since $\mu$ and $\nu$ have finite second moment, this is easily seen to be equivalent to $H(\mu |\gamma_d)<+\infty$ and $H(\nu |\gamma_d)<+\infty$, which is true according to Lemma \ref{finite-entropy}. 
According to Theorem \ref{thm:existence}, there exists a coupling $\pi^\ep (dxdy)= f^\ep(x)g^\ep(y)\,R^\ep(dxdy) \in C(\mu,\nu)$ such that $f^\ep$ is log-convex and $g^\ep$ is log-concave. It remains to show that this coupling is optimal for $\mathcal{T}^\ep_H(\mu,\nu)$. Since $P^\ep f^\ep(y)g^\ep(y) = e^{-W(y)}$, $y\in \R^d$, one sees that $\log g^\ep(y) = -W(y) - \log P^\ep f^\ep(y)$. The function $ \log P^\ep f^\ep$ is bounded continuous on the support of $\nu$ and $W$ is integrable with respect to $\nu$. Therefore, $\log g^\ep$ is integrable with respect to $\nu$. 
On the other hand, since $\log f^\ep \geq0$, the integral $\int \log f^\ep \,d\mu$ makes sense in $[0,+\infty]$. Let $\pi \in C(\mu,\nu)$ be a coupling such that $H(\pi|R^\ep)<+\infty$ (this set is non empty, since it contains $\mu\otimes \nu$). Applying Inequality \eqref{eq:dual} with $\alpha = \pi$, $\beta = R^\ep$ and $h(x,y)=\log f^\ep(x) + \log g^\ep(y)$, $x,y \in \R^d$, gives
\[
H(\pi | R^\ep) \geq \int \log f^\ep(x)+\log g^\ep(y)\,\pi(dxdy) = \int \log f^\ep(x)\,\mu(dx) + \int \log g^\ep(y)\,\nu(dy),
\]
which shows that $\log f^\ep$ is integrable with respect to $\mu$. A simple calculation shows that 
\[
H(\pi^\ep |R^\ep) = \int \log f^\ep(x)\,\mu(dx) + \int \log g^\ep(y)\,\nu(dy),
\]
which shows its optimality.
\endproof

Finally let us prove the technical lemmas used in the proof of Theorem \ref{thmCCT}.
\proof[Proof of Lemma \ref{finite-entropy}]
Let us first show that the probability measure $\mu$ has finite entropy. Since $\mu$ has a finite second moment, it is enough to show that $H(\nu | \gamma_d)<+\infty$, which amounts to show that $V$ is $\mu$ integrable. Since $V$ is bounded from below by some affine function, it is clear that $[V]_-$ is $\mu$ integrable. Moreover, since the convex function $V$ is such that $\int e^{V(x)}\,\gamma_d(dx) = 1$, this implies according to \cite[Lemma 2.1]{GMRS17} that $[V]_+(x) \leq \frac{|x|^2}{2}$, for all $x \in \R^d$, and so $[V]_+$ is also $\mu$ integrable. Similarly, to see that $H(\nu |\gamma_d)<+\infty$, it is enough to show that $W$ is $\nu$ integrable. On the one hand, $\int [W]_+\,d\nu = \int \left[\log(e^{-W})e^{-W}\right]_- \,d\gamma_d \leq \frac{1}{e}.$ On the other hand, $[W]_-$ is $\nu$ integrable since $W$ is bounded from below by some affine function.
\endproof

The following Lemma was used in the proof of Theorem 1: 

\begin{lem}\label{lem:technique}
Let $\nu(dx) = e^{-W(x)}\,\gamma_d(dx)$ with $W:\R^d \to \R\cup\{+\infty\}$ convex and $\eta \leq_c\nu$.  Assume furthermore that $\nu$ has compact support. 
Define, for all $\theta \in (0,\pi/2)$,  
\[
\nu_\theta = \mathrm{Law} (\cos \theta X + \sin \theta Z)\qquad \text{and} \qquad \eta_\theta = \mathrm{Law} (\cos \theta Y + \sin \theta Z),
\]
where $X \sim \nu$, $Y \sim \eta$ and $Z$ independent of $X, Y$ and such that the law $\alpha$ of $Z$ is given by $\alpha(dz) = \frac{1}{C} \mathbf{1}_{B}\,\gamma_d(dz)$ where $B$ is the Euclidean unit ball and $C$ a normalizing constant.
Then, for all $\theta \in (0,\pi/2)$,
\begin{enumerate}
\item the probability $\nu_\theta$ has a density of the form $e^{-W_\theta}$ with respect to $\gamma_d$, with $W_\theta : \R^d \to \R \cup\{+\infty\}$ convex,
\item the probability measures $\nu_\theta$ and $\eta_\theta$ are compactly supported,
\item it holds $\eta_\theta \leq \nu_\theta$,
\item the probability $\eta_\theta$ has finite entropy.
\end{enumerate}
\end{lem}
\proof
The density $f_\theta$ of $\nu_\theta$ is given by
\[
f_\theta(x) = \frac{1}{C'} \int_B e^{-W\left(\frac{x-\sin \theta y}{\cos \theta}\right)} e^{-\frac{|x-\sin \theta y|^2}{2\cos ^2\theta}} e^{-\frac{|y|^2}{2}}\,dy,
\]
where $C'$ is a normalizing constant. A simple calculation shows that
\[
e^{\frac{|x|^2}{2}}f_\theta(x) = \frac{1}{C'} \int_B e^{-W\left(\frac{x-\sin \theta y}{\cos \theta}\right)} e^{-\frac{|\sin \theta x - y|^2}{2\cos^2\theta}}\,dy
\]
and, according to Prekopa Theorem \cite{Pre73}, the right hand side is log-concave, which completes the proof of Item (1). The proofs of Items (2) and (3) are straightforward and left to the reader. 
The density of $\eta_\theta$ is $g_\theta (x) = \frac{1}{C}\int e^{-\frac{|\cos \theta y -x|^2}{2\sin^2 \theta} }\mathbf{1}_B\left(\frac{\cos \theta y -x}{\sin \theta}\right)\,\eta(dy)$ and so $g_\theta \leq \frac{1}{C}$. On the other hand, $g_\theta \log g_\theta \geq -1/e$. Since the support of $\eta_\theta$ is compact, one sees that $g_\theta \log g_\theta$ is integrable and so $\eta_\theta$ has finite entropy.
\endproof

\underline{\textbf{Acknowledgments}}: MF and NG were supported by ANR-11-LABX-0040-CIMI within the program ANR-11-IDEX-
0002-02. MF and MP were supported by Projects MESA (ANR-18-CE40-006) and EFI (ANR-17-CE40-0030) of the French National Research
Agency (ANR). We thank Christian L\'eonard for pointing out to us Fortet's work and sharing his notes explaining it, Gabriel Peyr\'e for his lectures on numerical optimal transport in March 2019 in Toulouse, and Franck Barthe and Michel Ledoux for useful discussions.

\bibliographystyle{amsplain}
\bibliography{bib}

\end{document}